\documentclass[12pt]{article}

\usepackage{epsfig}
\usepackage{amsmath}

\newcommand{\nit}{\noindent}
\newcommand{\no}{\nonumber}
\newcommand{\be}{\begin{equation}}
\newcommand{\ee}{\end{equation}}
\newcommand{\br}{\begin{eqnarray}}
\newcommand{\er}{\end{eqnarray}}

\newcommand{\ds}{\mbox{$\partial_{s}$}}

\newcommand{\lam}{\mbox{$\lambda$}}

\newcommand{\abs}[1]{\lvert #1 \rvert}
\newcommand{\norm}[1]{\lVert #1 \rVert}

\newtheorem{theo}{Theorem}[section]

\newtheorem{prop}{Proposition}[section]
\newtheorem{lem}{Lemma}[section]

\newtheorem{rmk}{Remark}[section]

\begin{document}

\title{Existence of KPP Type Fronts in Space-Time Periodic Shear Flows and 
a Study of Minimal Speeds Based on Variational Principle}

\author{James Nolen\thanks{Department of Mathematics, University of Texas at Austin, 
Austin, TX 78712 (jnolen@math.utexas.edu).}
\and Jack Xin\thanks{Department of Mathematics and ICES (Institute of Computational
Engineering and Sciences), University of Texas at Austin, Austin, TX 78712
(jxin@math.utexas.edu).}}

\date{}
\setcounter{page}{1}
\setcounter{section}{0}
\maketitle
\begin{abstract}
We prove the existence of reaction-diffusion traveling fronts in 
mean zero space-time periodic shear flows for
nonnegative reactions including the 
classical KPP (Kolmogorov-Petrovsky-Piskunov) nonlinearity. 
For the KPP nonlinearity, the minimal front speed is characterized 
by a variational principle involving the principal eigenvalue of 
a space-time periodic parabolic operator. 
Analysis of the variational principle shows 
that adding a mean-zero space time periodic shear flow to 
an existing mean zero space-periodic shear flow leads to 
speed enhancement. Computation of KPP minimal speeds is performed 
based on the variational principle 
and a spectrally accurate discretization  
of the principal eigenvalue problem. It shows that the enhancement 
is monotone decreasing in temporal shear frequency, and that 
the total enhancement from pure reaction-diffusion 
obeys quadratic and linear laws at small and large 
shear amplitudes.
\end{abstract}

\newpage

\section{Introduction}
\setcounter{equation}{0}
Front propagation in heterogeneous fluid flows is a fundamental issue 
in combustion science, and has been an active research area for decades 
(see \cite{CW}, \cite{Kh}, \cite{KA}, \cite{Ro},  
\cite{Vlad}, \cite{Xin1}, \cite{Yak} and references therein). 
A fascinating phenomenon is that the large time front speed 
depends strongly on the flow structures.  
Front dyanmics and speed enhancement have been studied mathematically 
for various flow patterns by analysis of the proto-type models, i.e. 
the reaction-diffusion-advection equations (see 
\cite{BH1,Const1,HPS,KR,Kh,MS1,MS2,NX1,PX,Vlad,Xin1,Xin2} and references 
therein). In \cite{Kh}, \cite{NX1}, fronts in 
space-time periodic shear flows (defined below) have been studied 
both analytically and numerically. This paper is a continuation of 
\cite{NX1} to establish front existence and variational principle 
on front speeds. The variational principle will be utilized both for 
analysis and computation to gain qualitative and quantitative understanding 
of front speeds. 

Let us consider the reaction-diffusion-advection equation:
\be
u_t  =  \Delta_x u + B\cdot \nabla_x u + f(u), \label{e1}
\ee
where $u = u(x,t)$, $x \in R^n$, $n \geq 2$, 
$t \in R$; $\Delta_x $ and $\nabla_x$ are the 
standard Lapliacian and gradient operators in $x$. 
To specify the vector field $B$ in $R^n$, let us write 
$x = (x_1 , y )$, with $ y \in  R^{n-1}$. Then  
$B$ is a smooth, time-dependent shear flow 
$B =  (b(y, t),0, \dots , 0)$ where $b(y,t)$ is a scalar function with 
period equal to one in both $y$ and $t$, and mean equal to zero over the 
period cell (denoted by $T^n$):
\be
\int_{T^n}\, b(y,t)\, dy \, dt =0. \label{e2}
\ee
To exclude degeneracy, we assume that:
\be 
\int_{T^n}\,|\nabla_y b(y,t)|^2 \, dy \, dt \not = 0. \label{e2a}
\ee 

We shall consider nonlinearity $f:[0,1] \to R$ to be smooth and of KPP type:
\be 
f(0) = f(1) = 0, \; f(u) > 0,\; \forall \; u \in (0,1), \; f'(1) < 0. 
\label{e3}  
\ee

\nit The definition (\ref{e3}) includes functions such as
$f(u) = u^m(1-u)$, $m \geq 1$, integer; and the Arrhenius-type 
nonlinearity $f(u) = e^{-E/u}(1-u)$, $E > 0$, arising in combustion 
\cite{Xin1}. The classical (strict) KPP nonlinearity \cite{Xin1} 
obeys the additional inequality $f(u) \leq u f'(0)$, $\forall \, u \in(0,1)$, 
e.g. $f(u)=u(1-u)$.

We are interested in traveling fronts 
moving in the direction along $x_1$ and of the form:
\be
u = U(x_1 - ct, y, t) \equiv U(s,y,\tau), \;\; s = x_1 - ct \in R,\; \tau = t, 
\label{e4}
\ee
with $U$ being 1-periodic in both $y$ and $\tau$ and 
satisfying conditions at $s = \pm \infty$:
\be
\quad \lim_{s \to -\infty} U(s,y,\tau) = 0, \;\;
 \lim_{s \to \infty} U(s,y,\tau) = 1, \label{e5}
\ee
uniformly in $y$ and $\tau$. After substitution, the problem becomes:
\br
\Delta_{s,y}\, U + (c + b(y,\tau))\, U_s - U_{\tau} + f(U) = 0,  \label{I2} \\
\lim_{s \to -\infty} U(s,y,\tau) = 0,  \; 
\lim_{s \to \infty} U(s,y,\tau) = 1, \; 
U \; {\rm has} \; {\rm period}\;  1\; {\rm in}\; (y,\tau), \no
\er
which has a parabolic equation on $R \times T^n$ with both the function 
$U$ and the constant $c$ unknown, or a nonlinear eigenvalue problem.  

The approach to front existence is to start with the front solutions to 
(\ref{e1}) with $f(u)$ being a so called 
combustion nonlinearity with cut-off: 
$ f_{\theta}(u) = 0 $ for $u \in [ 0, \theta ]$, 
$\theta \in (0,1)$, $f_{\theta}(u) > 0 $ for $u \in (\theta,1)$.
$f_{\theta}(1) = 0 $ and $f_{\theta}{'}(1) < 0$. 
In \cite{NX1}, the present authors outlined the 
existence proof of such combustion fronts using the method of 
continuation \cite{Xin3}. 
Here we provide the details, with the front speed estimates 
drawn from \cite{NX1}. The existence of 
KPP type fronts traveling at minimal speeds 
follows as suitable limits of the combustion fronts as $\theta \to 0$. 
Such a construction method 
of KPP type fronts is due to \cite{BerN1} on 
fronts in spatially dependent shears. 

In case of strict KPP nonlinearity, we show that the minimal speed denoted 
by $c^*$ is given by the variational principle:
\be
c^* = - \inf_{\lambda>0} \frac{\mu(\lambda)}{\lambda},  \label{e6}
\ee
where $\mu=\mu(\lambda)$ is the principal eigenvalue of the 
periodic parabolic operator 
\[ L^\lambda = \Delta_y \cdot - \partial_\tau \cdot + 
[\lambda^2 + \lambda \, b(y,\tau) + f'(0)],\;\; (y,\tau)\in T^n. \] 

The variational principle (\ref{e6}) allows us to deduce enhancement 
properties of $c^*$ from knowledge of $\mu$. Moreover, (\ref{e6}) provides 
the most efficient way of computing $c^*$ without solving the original 
time dependent PDE (\ref{e1}). We shall first compute 
$\mu$ by a spectrally accurate method, then minimize. 
With more accuracy in much shorter time, the minimal KPP front speeds 
are computed. The quadratic and linear speed enhancement laws 
from the pure reaction-diffusion speed (without advection), as well as 
monotonicity of enhancement in shear temporal frequency are 
demonstrated, in agreement with \cite{NX1} based on other methods. 
In particular, both analysis and computation from the 
variational principle show that 
adding a mean zero space-time shear to a spatial mean 
zero shear alway generates additional enhancement.


The paper is organized as follows. In section 2, we present properties and 
existence proof of combustion fronts. In section 3, we prove existence 
of KPP type fronts with minimal speeds. In section 4, we prove the 
variational principle (\ref{e6}) and its qualitative properties. 
In section 5, we show computational results of KPP minimal speeds 
based on (\ref{e6}).

\section{Combustion Fronts and Properties}
\setcounter{equation}{0}
In this section, $f$ is the combustion type nonlinearity with smooth cut-off. 
We prove:

\begin{theo}[Existence and Uniqueness] 
Let $f$ be of combustion type. Then there exists a classical solution $(U,c)$ 
to the problem (\ref{I2}). Moreover, 
this solution $(U,c)$ is unique up to translation in $s$. That is, if $(U_1,c_1)$ and $(U_2,c_2)$ are both solutions, then 
$c_1 = c_2$ and $U_1(s,y,\tau) = U_2(s+\lambda,y,\tau)$ for some $\lambda \in R$.
\end{theo}

The proof will be preceded by a list of properties of front solutions. 

\subsection{Strong Maximum Principle and Comparison}

As in \cite{BerN1} and \cite{Xin3}, one of the most important tools 
in our analysis is the strong maximum principle. 
We denote by $L_c$ the linear part of the 
operator in (\ref{I2}). This operator is parabolic, so the usual parabolic maximum principle and its corollaries hold. 
Specifically, if $U$ attains its minimum at 
some finite point $(s_0, y_0, \tau_0)$, then $U$ is constant for $\tau \leq \tau_0$. The functions on which $L_c$ operates are periodic in 
the $\tau$ variable (time), so the constant extends to all 
$\tau \in R$. This gives us a strong maximum principle for 
$L_c$ that is analogous to the maximum principle for elliptic operators:

\begin{prop}[Strong Maximum Principle]
Suppose $U \in C^2$ and 
\be
L_c u = \Delta_{s,y} U + (c + b)U_s - U_{\tau} \leq 0.
\ee
Suppose also that $U$ attains its minimum at some finite point $(s_0, y_0, \tau_0)$. Then $U$ is constant for all $(s,y,\tau) \in R \times T^n$. If $H(s,y,\tau ) \leq 0$ and $L_c u + H u \leq 0$ for all $(s,y, \tau )$ and $u$ 
attains a nonpositive minimum at some finite point, then 
$u$ is constant in $R \times T^n$.  Furthermore, if $L_c u + H u \leq 0$ and $u \geq 0$ in $R \times T^n$, then if $u(s_0, y_0, \tau_0) = 0$ 
at some finite point, then $u$ must be constant in $R \times T^n$, even if $H \leq 0$ does not hold.
\end{prop} 

Since solutions to (\ref{I2}) satisfy $L_c U \leq -f(U) 
\leq 0$, it follows that $0 < U < 1$. 
By the method of sliding domain \cite{BerN1}, \cite{Xin3}, the 
strong maximum principle implies the following two Theorems. 

\begin{theo}[Monotonicity] 
Let the nonlinearity $f$ be combustion type. 
Then the solution $U$ of (\ref{I2}) is strictly monotone 
in $s$: $U_s > 0$ for all $(s,y,\tau) \in R \times T^n$.
\end{theo}

\begin{theo}[Comparison Theorem]
Let $f \in C^1([0,1])$ be a combustion type nonlinearity. Let $\bar \beta(y,\tau)$ and $\beta(y,\tau)$ be smooth and periodic with 
$\bar \beta \geq \beta$ for all $(y,\tau) \in T^n$. 
Assume that $u$ and $\bar u$ $\in C^{2,1}_{(s,y),\tau}(R^1\times T^n)$ 
are two functions that satisfy
\br
L u + f(u) & = & \Delta_{s,y} u + \beta(y,\tau)u_s - u_\tau + f(u) \geq 0,  \no \\
\bar L \bar u + f(\bar u) & = &\Delta_{s,y} \bar u  + \bar \beta(y,\tau)\bar u_s - \bar u_\tau + f (\bar u) \leq 0, \label{SM2}
\er
and the conditions at $s$ infinities:
\br
\lim_{s \to -\infty} \bar u &=& 0 \;\;\; \text{and} \;\;\; \lim_{s \to +\infty} \bar u = 1  \no \\
\lim_{s \to -\infty} u &=& 0 \;\;\; \text{and} \;\;\; \lim_{s \to +\infty} u = 1 \no
\er
Additionally, assume that at least one of the following holds:
\be
\bar \beta \equiv \beta \quad \text{OR} \quad \bar u_s > 0 \quad \text{OR} \quad u_s > 0.
\ee
Then there exists $\lambda \in R$ such $u(s,y,\tau) = \bar u(s + \lambda,y,\tau)$ for all $(s,y,\tau)$. Furthermore, there exists a point $(y_0,\tau_0) \in T^n$ such that $ \bar \beta(y_0,\tau_0) = \beta(y_0,\tau_0)$, so that $\bar \beta > \beta$ cannot hold for all $(y,\tau) \in T^n$. 
\end{theo}

\nit {\it Proof:} 
The proof of Theorem 2.3 uses the sliding domain method. Define the function  
$w_\lambda = \bar u(s+\lambda,y,\tau) - u(s,y,\tau)$.
The key is that the inequalities (\ref{SM2}) are invariant under translation in the $s$ direction. For all $\lambda \in R$,
\be
L w_\lambda + H_\lambda w_\lambda \leq (\beta - \bar \beta)\bar u_s(s+\lambda) \label{SM3}
\ee
where
\be
H_\lambda = H_\lambda(s,y,\tau) = \int_0^1 f'(\xi \bar u(s+\lambda) + (1-\xi) u)\;d\xi.
\ee
Assuming either $\bar \beta \equiv \beta$ or $\bar u_s > 0$ is satisfied, the right hand side of (\ref{SM3}) is non-positive. In the event that $u_s > 0$, we consider instead
\be
\bar L w_\lambda + H_\lambda w_\lambda \leq  (\beta - \bar \beta) u_s(s+\lambda),
\ee
so that again the right side is non-positive. In either case, $w_\lambda$ satisfies
\be
\lim_{s \to -\infty} w_\lambda = 0 \;\;\; \text{and} \;\;\; \lim_{s \to +\infty} w_\lambda = 0. \no \\
\ee
Using the strong maximum principle for $L$ (or $\bar L$), we can show that for $\lambda$ sufficiently large, $w_\lambda > 0 $ everywhere in $R \times T^n$. Then, for $\mu$ defined by
\be
\mu = \inf \left\{ \lambda_0 \;\lvert \; w_\lambda > 0 \;\;\text{in}\;\; R \times T^n ,\;\;\forall \; \lambda > \lambda_0 \;\right\} , \no
\ee 
it can be shown that $w_\mu \equiv 0$. In other words, we can translate $\bar u$ to the left so that $\bar u > u$. Then translating back to the right until they touch, 
we have $\bar u \geq u$ with equality holding at a finite point, and the strong maximum principle implies the result. If $\bar \beta > \beta$, the right hand side of 
inequality (\ref{SM3}) becomes strictly less than zero for all $\lambda$. Plugging in $w_\mu$, we then obtain a contradiction, implying $\bar \beta(y_0,\tau_0) = \beta(y_0,\tau_0)$ at 
some point. The proof relies on the fact that $f(u)$ vanishes for $u< \theta$ and that $f'(1) < 0$, so that $H_\lambda \leq 0$ in the region where $\abs{s}$ is sufficiently large. 
Note that this is not necessarily true for the KPP nonlinearity because in that case $f(u) > 0$ for all $u \in (0,1)$.

Theorem 2.2 follows from Theorem 2.3. We take $\bar u = u$ and 
$\bar \beta = \beta$ and show that $\mu = 0$, hence, $u_s \geq 0$. 
Differentiating equation (\ref{I2}) with respect to $s$, 
we obtain an equation for $u_s$. The strong maximum principle and 
the limits as $s \to \pm \infty$ then imply strict inequality: $u_s > 0$.

\subsection{Exponential Decay}
As with the problems considered in \cite{Xin3} and \cite{BerN1}, we can obtain precise information about the decay rate of solutions through related eigenvalue problems. Specifically, 
solutions of (\ref{I2}) must decay exponentially fast to the limits 0 and 1 as $s \to \pm \infty$, as described in the following theorem:
\begin{theo}[Exponential Decay]
If (U,c) is a classical solutions of (\ref{I2}) with combustion nonlinearity, then there exist constants, $C_1 >0$, $C_2 > 0$, $\lam_1 >0$, $\lam_2 <0$,  and positive functions 
$\Phi_1 (y,\tau), \; \Phi_2 (y,\tau) \in C^2(T^n)$ such that
\br
U(s,y,\tau) & \leq & C_1 e^{\lam_1 s}\Phi_1 (y,\tau), \quad s \leq s_1  \no \\
1 - U(s,y,\tau) & \leq & C_2 e^{\lam_2 s}\Phi_2 (y,\tau), \quad s \geq s_2
\er
\end{theo}
\nit {\it Proof:}
This result is analogous to Proposition 1.1 in \cite{Xin3}, the only significant difference being that the related eigenvalue problem 
is parabolic. If $s_1 < 0$ is such that $U(s,y,\tau) < \theta$ for $s \leq s_1$, then we have
\be
L_c U = 0 \quad \text{for} \; \; s \leq s_1 .
\ee
To compare $U$ with an exponentially decaying function, we want to solve $L_c w = 0$ for some function $w$ of the 
form $w = C_1 e^{\lam s}\Phi(y,\tau)$ and $\Phi > 0$.  Such a function must satisfy
\be
e^{-\lam s}L_c w= \Delta_y \Phi - \Phi_{\tau} + [\lam^2 + \lam(c + b(y,\tau))]\Phi = 0.
\ee

We will call $\rho(\lambda) \in R$ the principal eigenvalue of this operator if for some $\Phi > 0$, 
\be
\Delta_y \Phi - \Phi_{\tau} + [\lam^2 + \lam(c + b)]\Phi = \rho(\lam)\Phi . \label{EX1}
\ee
Thus, we want to find some $\lam >0$ such that $\rho(\lam) =0$.

\begin{lem}
For each $\lambda \in R$ there exists a principal eigenvalue $\rho(\lambda)$ and a principal eigenfunction $\Phi > 0$ solving (\ref{EX1}).
\end{lem}

\nit {\it Proof:}
This is a periodic-parabolic eigenvalue problem, and again we see the importance of the periodicity of the time dependence that leads to a strong maximum principle. 
Problems of this type have been treated in more detail in \cite{Hess} and \cite{Laz}. Let 
\be
E = C^{2+\alpha, 1+\alpha/2}(T^n) \;\;\; \text{and} \;\;\; F = C^{\alpha, \alpha/2}(T^n) .
\ee
Let $M:E \to F$ be the operator defined by, 
\be
M \Phi =  \Phi_{\tau} -\Delta_y \Phi  - [\lam^2 + \lam(c + b)]\Phi .
\ee
For fixed $\lambda$ we can pick $d >0$ sufficiently large so that the operator $(M+d) \Phi $ is invertible (see Section 2 of \cite{Laz}).  
By parabolic Schauder estimates, the map $T = i(M+d)^{-1}:F \to F$ is compact, where $i$ is the embedding of 
$E$ into $F$.  Since the strong maximum principle implies positivity of the operator $T$, the Krein-Rutman theorem then 
implies the existence of a simple eigenvalue $\mu$ with positive 
eigenfunction $\Phi$ for operator $T$. 
Hence $T \Phi = \mu \Phi$, or $M \Phi = (\frac{1}{\mu} - d) \Phi $, so that $\rho(\lambda) = \frac{1}{\mu} - d$ and $\Phi$ are the 
desired eigenvalue and eigenfunction for $M$.

Because of the periodicity in $y$ and $\tau$, integrating the equation over $T^n$ gives us
\be
\lam^2 \int_{T^n} \Phi \,dy\,d\tau + \lam \int_{T^n} (c + b)\Phi \,dy\,d\tau = \rho(\lam) \int_{T^n} \Phi \,dy\,d\tau , \label{EX2}
\ee
and since $\Phi > 0$ we can divide by $\int_{T^n} \Phi$ to obtain
\be
\rho(\lambda) = \lam^2  + \lam c +  \lambda \frac{\int_{T^n} b\Phi}{\int_{T^n}\Phi} . \label{EX3}
\ee
The quotient on the right is bounded in absolute value by $\lambda \norm{b}_\infty$, so we see that as $\lam \to \infty$, $\rho(\lam) \to +\infty$. 
Setting $\lam = 0$ it is also clear that we must have $\rho(0)=0 $ 
with $\Phi \lvert_{\lambda = 0} \equiv 1$ (i.e. positive constant).  Taking the derivative of equation (\ref{EX1}) with respect to $\lam$ and 
setting $\lam = 0$ and $\Phi \lvert_{\lambda = 0} \equiv 1$ we obtain
\be
\Delta_{y} \Phi_{\lam} - (\Phi_{\lam})_{\tau} + (c + b) = \rho_{\lam}(0) .
\ee
Averaging over $y$ and $\tau$ we get
\be
\rho_{\lam}(0) = c ,
\ee
since  $b(y,\tau)$ is mean zero and  $\Phi$ is periodic.  This implies that
\be
\rho(\lam) = c\lam + O(\lam^2) .
\ee
If we integrate equation (\ref{I2}) over all $s,y,\tau$, then we see because of periodicity, integration by parts gives us
\be
c = -\int_{R \times T^n} f(U) < 0 .
\ee
Thus, for $\lam$ small and positive, we must have $\rho(\lam) <0 $.  Since $\rho(\lam) \to +\infty$ as $\lam \to \infty$, there must be some 
$\lam_1 >0$ such that $\rho(\lam_1) = 0$.  The principal eigenfunction, $\Phi_1$, corresponding to eigenvalue $\rho(\lam_1)$ is strictly positive, and
$w = e^{\lam_1 s}\Phi_1$ solves $ \bar L w = 0$.

To prove exponential convergence as $s \to +\infty$, we perform the same analysis on the function $V = 1-U$ in the region $[s_2,+\infty) \times T^n$, where $s_2>0$ is such that
\be
f'(U) \leq \frac{1}{2} f'(1) < 0, \; \; \; \forall \; s>s_2 .
\ee

\subsection{Existence of Traveling Front Solutions}
Now that we have have established the strong maximum principle for $L_c$, the exponential behaviour of solutions as $s \to \pm \infty$, and the existence of eigenvalues for the related eigenvalue problem (\ref{EX1}),  the proof of Theorem 2.1 follows from the reasoning in \cite{Xin3} with modification to 
account for the time-dependent component of operator. We work with the equation in the whole space $R \times T^n$ and prove existence of solutions via the method of continuation.

First, supposing that for some $b = b_0$ we have a solution $(U_0,c_0)$ to the problem (\ref{I2}), we find that there exist solutions corresponding to small perturbations:
\be
b^{\delta}  =  b_0 + \delta b_1^{\delta} . \no
\ee
We look for solutions in the same functions spaces as in \cite{Xin3} with the addition of the $\tau$ dimension to the domain of the functions. As in \cite{Xin3}, we consider the linearized operator $L_{c_0} v + f'(U_0)v$ and find that, because of monotonicity of solutions to (\ref{I2}), this linearized operator is Fredholm with one-dimensional kernel. Then, the existence of a solution $(U^{\delta}, c^{\delta})$ corresponding to shear $b^\delta$ follows from a contraction mapping principle for $\delta$ small.

Next, to complete the continuation, we suppose there is some maximum value $D>0$ such that we can find solutions for all $\delta \in [0,D)$, and we want to obtain a solution for $\delta = D$, 
thus implying that no such maximum exists. To do so, we consider a sequence of solutions $(U^{\delta}, c^{\delta})$ to the problems:
\br
\Delta_{s,y} U^{\delta} + (c^{\delta} + b^{\delta})U^{\delta}_s - U^{\delta}_{\tau} + f(U^{\delta}) = 0 ,  \label{Conv1} \\
\lim_{s\to\infty} U^\delta = 1,\; 
\lim_{s\to -\infty} U^\delta = 0,\; 
U^\delta \; {\rm has} \; {\rm period}\;  1\; {\rm in}\; (y,\tau). \label{I4.1} 
\er
We show that there is a subsequence $(U^{\delta_j},c^{\delta_j})$ converging to a solution in the case $\delta = D$.  First, we need to bound the wave speeds $c^\delta$ in 
order to extract a convergent subsequence. 
To do this, we compare these speeds to the wavespeed of the corresponding one-dimensional problem:
\br
\ddot v + k \dot v + f (v) = 0, \label{Conv2} \\
\lim_{s\to\infty} v = 1,\; \lim_{s\to -\infty} v = 0, \no
\er
with normalization
\be
v(0) = \frac{1}{2} .
\ee

\begin{lem}
Let $(U^\delta, c^\delta)$ solve (\ref{Conv1}) and let $(v,k)$ solve (\ref{Conv2}).  Then $k \leq c_\delta + \norm{b^\delta}_\infty$.
\end{lem}

\nit {\it Proof:}
Defining $\bar u(s,y,\tau) = v(s)$ and $u = U^\delta$, we have 
\br
\Delta_{s,y} u + \beta u_s - u_\tau + f(u) = 0, \no \\
\Delta_{s,y} \bar u + \bar \beta \bar u_s - \bar u_\tau + f(\bar u) = 0, \no
\er
with $\bar \beta(y,\tau) \equiv k$ and $\beta = c^\delta + b^\delta (y,\tau)$. If $k \geq c_\delta + \norm{b^\delta}_\infty$, then $\bar \beta \geq \beta$, so we can apply Theorem 2.3 to conclude that there must be some $(y_0,\tau_0) \in T^n$ such that $\bar \beta(y_0,\tau) = \beta(y_0,\tau_0)$.  This implies $k = c^\delta + b^\delta(y_0,\tau_0) \leq c^\delta + \norm{b^\delta}_\infty$.

Now for $\delta \in [0,D)$, Lemma 2.2 implies there exists constant $K<0$ independent of $\delta$ such that $K < c^{\delta} < 0$.  Hence, 
we can pick a subsequence $\{ c^{\delta_n} \}$ that converges to some $c^* \leq 0$. To have a solution $(U,c)$ with $\delta = D$, however, we must have $c^* < 0$. The following lemma was proved in the appendix of \cite{NX1}:

\begin{lem}[Bounding $c_*$ away from zero] 
If $\{ (U^{\delta},c^{\delta}) \}$ is a set of solutions to (\ref{Conv1}) with $\delta \in [0,D)$, then there is a constant $\gamma > 0$ independent of 
$\delta $ such that
\be
c^{\delta} < - \gamma < 0, \;\;\;\; \forall \;\; \delta \in [0,D).
\ee
Hence, $c^* < 0$.
\end{lem}

Now, by standard interior parabolic estimates (\cite{Lad}, \cite{Lbm}), we can find a function $U^*$ and a subsequence $\{ U^{\delta_j} \}$ such that 
on any compact set $\Omega \subset R \times T^n$, $U ^{\delta_j} \to U^*$ uniformly with $U^*$ solving
\be
\Delta_{s,y} U^* + (c^* + b^{D})U^*_s - U^*_{\tau} + f(U^*) = 0 . \label{Conv3}
\ee

We need to show that the limits (\ref{I4.1}) hold for $U^*$. Since solutions are invariant with respect to translation in the $s$ variable (that is, we can translate them and they are still solutions), we can 
impose the normalization condition
\be
\max_{(y,\tau \in T^n)} U^{\delta}(0,y,\tau) = \theta
\ee
on all solutions $(U^{\delta},c^{\delta})$.  By monotonicity in $s$, $U^{\delta}(s,y,\tau) \leq \theta$ for $(s,y,\tau) \in (-\infty,0) \times T^n$.  Moreover, 
the limiting function $U^*$ must satisfy $U^*_s \geq 0$ for all $(s,y,\tau)$.

First, we show that $\lim_{s \to -\infty} U^*(s,y,\tau) = 0 $. In the proof of Theorem 2.4, we showed that for each $\delta$ we have some $C^{\delta},\lam^{\delta}, \Phi^{\delta}$ such that 
\be
U^{\delta}(s,y,\tau) \leq  C^{\delta} e^{(\lam^{\delta} s)}\Phi^{\delta} (y,\tau), \quad s \leq 0  \no \\
\ee
(here we take $s_1$ to be 0 because of the normalization condition). Because the eigenfunctions $\Phi^\delta$ are unique up to multiplication by a constant, we may assume that 
$\int_{T^n} \Phi^\delta = 1$ for each $\delta$. Since the wave speeds are bounded away from zero, $c^{\delta} < - \gamma < 0$, 
the corresponding $\lambda^\delta$ must be bounded away from zero. Supposing this were not the case, we take a sequence $\left\{ \delta_j \right\}^{\infty}_{j=1}$ such that $\lambda^{\delta_j} \to 0$ and 
$\Phi^{\delta_j} \to \Phi_0$ in $C^1(T^n)$ with $\Phi_0$ solving 
\be
\Delta_{y} \Phi - \Phi_{\tau} = 0 .
\ee
The strong maximum principle implies that $\Phi_0 \equiv 1$, which means that $\int_{T^n} b \Phi_0 = 0$ since $b$ is mean zero. However, using (\ref{EX3}) and $\int_{T^n} \Phi^{\delta_j} = 1$, we have
\be
\lam^{\delta_j} = -c^{\delta_j} - \int_{T^n} b \Phi^{\delta_j} \geq \gamma - \int_{T^n} b \Phi^{\delta_j}.
\ee
Letting $j \to \infty$ we obtain a contradiction since the integral on the right hand side goes to zero. Thus, there exists a constant $\lambda^* > 0$ such that 
$\lambda^\delta > \lambda^*$ for all $\delta$. Harnack's inequality implies that there is a constant $M > 0$ independent of $\delta$ such that
\be
\inf_{(y,\tau)} \Phi^{\delta} \geq M \sup_{(y,\tau)} \Phi^\delta \geq M .
\ee
Now chosing $C^* = 1/M$, we have
\be
U^{\delta}(s,y,\tau) \leq C^* e^{\lam^* s} \;\; \text{whenever}\;\;\; s<0,\;\forall \;\;(y,\tau) .
\ee
This implies the desired result: $\lim_{s \to -\infty} U^*(s,y,\tau) = 0 $.

Now we show that $\lim_{s \to \infty} U^*(s,y,\tau) = 1 $.  Since $U^*_s \geq 0 $, we can define
\be
\lim_{s \to +\infty} U^*(s,y,\tau) = \psi(y,\tau) ,
\ee
and we want $\psi \equiv 1$.  It is not hard to see that the function $\psi$ must solve the equation
\be
\Delta_{y} \psi - \psi_\tau + f(\psi) = 0
\ee
for $(y,\tau) \in T^n$. To see this, let $h(y,\tau)$ be a smooth function on $T^n$. Also, let $\xi(s):R \to [0,1]$ be smooth with $\xi \equiv 0$ for $s < 1$ and $\xi \equiv 1$ for $s > 2$. 
  Multiply equation (\ref{Conv3}) by the test function $ \phi(s,y,\tau) = h(y,\tau) \xi(s) $, average over $[0,L] \times T^n$: 
\be
\frac{1}{L} \int_0^L \int_{T^n} \Delta_{s,y} U^* \phi + (c + b)U^*_s \phi + U^*_\tau \phi + f(U) \phi = 0 . \label{LIM1}
\ee
Parabolic estimates imply the boundedness of the $s$ derivative of $U^*$. As a result, integrating by parts and taking the limit as $L \to \infty$ gives us the result
\be
\int_{ T^n}  \psi \Delta_{y} h  + \psi h_\tau + f(\psi) h = 0
\ee
for any test function $h \in C^{\infty}(T^n)$. Thus, $\psi$ solves 
\be
\Delta_{y} \psi - \psi_\tau + f(\psi) = 0
\ee
for $(y,\tau) \in T^n$. Now, the maximum principle implies that $\psi$ must be constant.  Hence,
\be
\Delta_y \psi - \psi_\tau \equiv 0 = -f(\psi),
\ee
so $\psi \equiv \theta $ or $\psi \equiv 1$, since $f(\psi) > 0$ for $ \psi \in (\theta, 1)$.  If $\psi \equiv \theta $, then by monotonicity, $U^* \leq \theta $ for all 
$(s,y,\tau)$ so that $U^*$ satisfies
\be
\bar L u = \Delta_{s,y} U + (c + b^{D})U_s - U_{\tau} = - f(U) \equiv 0 .
\ee
Also, $U^*$ attains its maximum at some finite point $(0,y_0,\tau_0)$ due to the normalization condition.  So, the strong maximum principle implies that
\be
U^* \equiv \theta ,
\ee
which contradicts the limiting behavior as $s \to - \infty$.  Therefore, $\psi \equiv 1$ and we have shown that 
\be
\lim_{s \to \infty} U^*(s,y,\tau) = 1 .
\ee

Now that we have a solution $(U,c)$ to problem (\ref{I2}), the uniqueness of $U$ (up to translation in $s$) follows easily from the comparison theorem. Let $(U,c)$ and $(\bar U,\bar c)$ be two solutions to (\ref{I2}) with $\bar c \geq c$. Theorem 2.3 implies that we must have $c = \bar c$ and that for some constant $\lambda \in R$, $U(s,y,\tau) = \bar U(s + \lambda,y,\tau)$ for all $(s,y,\tau) \in R \times T^n$.

\begin{rmk}
Although periodicity in time plays a crucial role in the proof of these results, the periodicity of solutions in the $y$ variable is not necessary. That is, if $\omega \subset R^{n-1}$ is a bounded region with smooth boundary, our results generalize easily to treat traveling fronts in the cylinder $(x,y) \in (R \times \omega)$ when the shear is periodic in time with Neumann boundary conditions imposed in the $y$ dimension: $\frac{\partial u}{\partial \nu} = 0$ for $y \in \partial \omega$. With the Hopf Lemma, the strong maximum principle of Proposition 2.1 still holds, as well as the rest of the results of this section. The same is true of our results for the KPP equation described in following sections.
\end{rmk}

\section{Existence of KPP Type Fronts}
\setcounter{equation}{0}
Using the results for the combustion case, 
we now construct a solution to the problem
\br
\Delta_{s,y} U + (c + b(y,\tau))U_s - U_\tau + f(U) = 0, \label{KP1} \\
\lim_{s\to\infty} U = 1, \; \lim_{s\to -\infty} U = 0, \; 
U \; {\rm has} \; {\rm period}\;  1\; {\rm in}\; (y,\tau), \label{KP1L} 
\er
where $f(u)$ is KPP type nonlinearity (\ref{e3}). 
Since the shear $b(y,\tau)$ will be the same throughout this section, 
we will use the notation $L_c $ to denote the operation
\be
L_c U = \Delta_{s,y} U + (c + b(y,\tau))U_s - U_\tau . \no
\ee 
 
\begin{theo}[Existence] 
Let $f \in C^1([0,1])$ with $f(u) > 0$ in $(0,1)$,  $f(0) = f(1) = 0$. Then there exists a classical solution $(U^*,c^*)$ to the 
problem (\ref{KP1})-(\ref{KP1L}). 
\end{theo} 
\nit {\it Proof:}
This proof extends the arguments in section 8 of \cite{BerN1} to the case with time-dependent coefficients. Following \cite{BerN1}, we let $\chi_\theta(u) \in [0,1]$ be a smooth cutoff function satisfying
\br
\chi_\theta(u) & = & 0 \quad \text{if} \;\; u \leq \theta/2 , \no \\
\chi_\theta(u) & = & 1 \quad \text{if} \;\; u \geq \theta , \no
\er
and $\chi_{\theta'} \geq \chi_\theta$ whenever $\theta' < \theta $. Each function $f_\theta(u) = f(u) \chi_\theta(u)$ is a nonlinearity of combustion 
type and $f_\theta \to f$ uniformly on $[0,1]$ as $\theta \to 0$. We will also assume that if $\theta' < \theta$, then for some point $u_0 \in (0,\theta)$, $\chi_{\theta'}(u_0) > \chi_{\theta}(u_0)$ and thus $f_{\theta'}(u_0) > f_{\theta}(u_0)$. By Theorem 2.1, for each $\theta \in (0,1/2)$ there is a unique pair $(u_\theta, c_\theta) $ solving
\br
L_{c_\theta} u_\theta + f_\theta (u_\theta) = 0, \label{KP2} \\
\lim_{s\to\infty} u_\theta = 1, \quad \lim_{s\to -\infty} u_\theta = 0, \no
\er
with normalization condition
\be
\max_{(y,\tau)} u_\theta (0,y,\tau) = \frac{1}{2} . \no
\ee
We will show that the solutions $(u_\theta , c_\theta)$ converge to a solution $(u^*,c^*)$ of the KPP problem and that this speed $c^*$ is the 
minimum possible speed (minimal in absolute value) of any solution to (\ref{KP1}). First, 
we bound the sequence of wave speeds in order to extract a convergent subsequence.

\begin{lem} 
The sequence of speeds $c_\theta$ is strictly decreasing as $\theta \to 0$. Moreover, the sequence is bounded, so that there exists $c^*<0$ such that $lim_{\theta \to 0}\, c_{\theta} = c^*$. 
\end{lem}
\nit {\it Proof:}
For the first claim, let $\theta_2 < \theta_1 < \frac{1}{2}$ with corresponding solutions $(u_1,c_1)$ and $(u_2,c_2)$:
\be
L_{c_1} u_1 + f_{\theta_1} (u_1) = 0, \no 
\ee
and
\br
L_{c_2} u_2 + f_{\theta_2} (u_2) &=& 0, \no \\
L_{c_2} u_2 + f_{\theta_1} (u_2) &\leq& 0.\no
\er
Arguing by contradiction, suppose $c_2 \geq c_1$. Then for $\bar \beta = c_2 + b(y,\tau)$ and $\beta = c_1 + b(y,\tau)$, we have $\bar \beta \geq \beta$, so we can apply Theorem 2.3 
to conclude that $\bar \beta(y_0,\tau_0) = \beta(y_0,\tau_0)$ for some point $(y_0,\tau_0)$ and that $u_1(s,y,\tau) \equiv u_2(s + \lambda,y,\tau)$ for some translation $\lambda \in R$. 
This implies that $c_1 = c_2$ and that $u_1$ must satisfy both
\br
L_{c_1} u_1 + f_{\theta_1} (u_1) &=& 0,  \\
L_{c_1} u_1 + f_{\theta_2} (u_1) &=& 0.
\er
This cannot be true, however, since for some value of $u_1$, $f_{\theta_1}(u_1) < f_{\theta_2}(u_1)$, by construction of the functions $f_{\theta_1}$,$f_{\theta_2}$. Therefore $c_2 < c_1$.

As in \cite{BerN1}, we obtain a lower bound on the sequence $\{ c_\theta \}$ by a comparison with the wavespeeds of the corresponding 1-D problems:
\br
\ddot v_\theta + k_\theta \dot v_\theta + f_\theta (v_\theta) = 0, \\
\lim_{s\to\infty} v_\theta = 1, \quad \lim_{s\to -\infty} v_\theta = 0, \no
\er
with normalization
\be
v_\theta(0) = \frac{1}{2} .
\ee
Using $k = k_\theta$ and $f = f_\theta$, it follows from Lemma 2.2 that for 
each $\theta \in (0,1/2)$, $k_\theta \leq c_\theta + \norm{b}_\infty$. Thus, a lower bound on $k_\theta$ implies a lower bound on $c_\theta$. Such a bound was derived in section 8 of \cite{BerN1}.

Since all the $f_\theta$ and $u_\theta$ are uniformly bounded, then as before, parabolic estimates imply the existence of a subsequence $u_{\theta_k}$ converging uniformly on compact sets to a function $U^*$ that solves (\ref{KP1}).

We must show, however, that the limits (\ref{KP1L}) are achieved. We know that for each $\theta$, $u_\theta$ is 
monotone in $s$. This implies $U^*_s \geq 0$. So we can define the limits
\be
\psi^+ (y,\tau) = \lim_{s \to \infty} U^*(s,y,\tau) \;\; \text{and} \;\; \psi^- (y,\tau) = \lim_{s \to -\infty} U^*(s,y,\tau) .
\ee
As in equation (\ref{LIM1}), we multiply by a smooth test function and integrate by parts to show that that $\psi^+$ and $\psi^-$ satisfy
\br
\Delta_{y} \psi^+ -\psi^+_\tau = -f(\psi^+),  \no \\
\Delta_{y} \psi^- -\psi^-_\tau = -f(\psi^-) .
\er
Integrating both equations over $T^n$ shows that $f(\psi^+) \equiv 0$ and $f(\psi^-) \equiv 0$. so that $\psi^{\pm} \equiv \text{const.}$, 
i.e. $0$ or $1$. Because of the normalization condition on the functions $U_\theta$, however, we must have 
\be
\max_{(y,\tau) \in T^n} U^*(0,y,\tau) = 1/2,
\ee 
so that $\psi^- \equiv 0$ and $\psi^+ \equiv 1$. Thus, the limits (\ref{KP1L}) as $s \to \pm \infty$ are acheived, concluding the proof of Theorem 3.1.

In the following sections, we will use $c^*$ to denote this unique speed obtained by the above construction of the solution $(U^*,c^*)$.

\begin{prop} 
If $(u,c)$ is any solutions to (\ref{KP1}) and (\ref{KP1L}), then $c \leq c^*<0$.
\end{prop}
\nit {\it Proof:}
Suppose $(\bar u,c)$ is a solution to (\ref{KP1}). For any $\epsilon > 0$, we can pick $\theta> 0$ sufficiently small such that $c_\theta \in (c^*,c^* + \epsilon)$, where $(u_\theta, c_\theta)$ is the 
corresponding solution to the problem with combustion nonlinearity (\ref{KP2}), as in the construction of $(u^*,c^*)$. Because $f_\theta \leq f$, we have
\br
& & L_c \bar u + f(\bar u) = 0 , \no \\
& & L_c \bar u + f_\theta(\bar u) \leq 0 ,  \no 
\er
and
\be
L_{c_\theta} u_\theta + f_\theta (u_\theta) = 0. \no
\ee
Arguing by contradiction, assume $c \geq c_\theta$. 
We know that $(u_\theta)_s > 0$, so we can apply Theorem 2.3 with $\bar \beta = c + b(y,\tau) \geq c_\theta + b(y, \tau) = \beta$ to conclude that 
$c = c_\theta$.  Hence $c \leq c_\theta < c^* + \epsilon$. 
Letting $\epsilon \to 0$, we conclude that $c \leq c^*$.

\section{Variational Characterization of $c^*$}
\setcounter{equation}{0}
The existence results in the previous section hold for a general KPP nonlinearity. If we impose the additional condition $f(u) \leq uf'(0)$ for all $u \in[0,1]$, we can obtain a variational characterization of the minimal wavespeed $c^*$ in terms of a related eigenvalue problem on $T^n$, extending the characterizations given in \cite{BerN1} and \cite{BH1} to this problem with a time-dependent shear.

Let $\mu(\lambda)$ be the principal eigenvalue of the operator $L^\lambda = \Delta_y \cdot - \delta_\tau \cdot + [\lambda^2 + \lambda b(y,\tau) + f'(0)]$ 
with associated principal eigenfunction $\phi(y,\tau) > 0$. That is,
\be
L^\lambda \phi = \Delta_y \phi + [\lambda^2 + \lambda b(y,\tau) + f'(0)]\phi - \phi_\tau = \mu(\lambda) \phi  \label{VAR1}
\ee
The variational characterization is given in the following theorem:
\begin{theo}
If the nonlinearity $f$ is KPP-type and satisfies $f(u) \leq uf'(0)$ 
for all $u \in [0,1]$, then the minimal speed $c^*\, < \, 0$ 
is given by the variational formula
\be
c^* = - \inf_{\lambda>0} \frac{\mu(\lambda)}{\lambda} .  \label{VAR2}
\ee
\end{theo}

\nit {\it Proof:}
We define the value 
\be
\gamma^* = - \inf_{\lambda>0} \frac{\mu(\lambda)}{\lambda} .\label{Var3}
\ee
Clearly, $\gamma^* > -\infty$. Let $c<0$ satisfy $c<\gamma*$. We show first that we must have $c \leq c^*$, 
implying that $\gamma^* \leq c^*$. Since $c < \gamma^*$, there exists $\lambda_0 > 0$ such that
\be
-\lambda_0 c > \mu(\lambda_0) . \no
\ee
Letting $\phi_0 = \phi_0(y,\tau) > 0$ be the principal eigenfunction associated with $\lambda_0$ we have
\be
\Delta_y \phi_0 + [\lambda_0^2 + \lambda_0 b(y,\tau) + f'(0)]\phi_0 - (\phi_0)_\tau = \mu(\lambda_0) \phi_0 < -\lambda c \phi_0 , \no
\ee
or 
\be
\Delta_y \phi_0 + [\lambda_0^2 + \lambda_0 (c + b(y,\tau)) + f'(0)]\phi_0 - (\phi_0)_\tau < 0 .
\ee
Define $\bar u(s,y,\tau) = e^{\lambda_0 s}\phi_0(y,\tau)$, and we see that $\bar u$ solves
\br
L_c \bar u + f'(0)\bar u &<& 0  \no \\
L_c \bar u + f(\bar u) &<& 0 . \label{VAR4}
\er

To show that $c \leq c^*$, we argue by contradiction and suppose $c^* < c$.  As in the construction of the solution $(u^*,c^*)$, pick $\theta > 0$ sufficiently 
small so that $c^* < c_\theta < c < \gamma^*$ where $(u_\theta,c_\theta)$ solves (\ref{KP2}).  Since $c_\theta < c $ and $f_\theta \leq f$,
\br
L_{c_\theta} u_\theta + f_\theta (u_\theta) = 0  \no \\
L_c \bar u + f_\theta (\bar u) < 0 . \label{VAR5}
\er
Replacing $\bar u$ with the function $\min(\bar u, 1)$, we can apply Theorem 2.3 to $u = u_\theta$ and $\bar u$ with $\beta = c_\theta + b(y,\tau)$ 
and $\bar \beta = c + b(y,\tau)$ to conclude that for some $(y_0,\tau_0)\in T^n$, $\bar \beta(y_0,\tau_0) = \beta(y_0,\tau_0)$, hence $c_\theta = c$, a contradiction. 
Therefore, $c \leq c^*$, which implies $\gamma^* \leq c^*$.

Now the opposite inequality $\gamma^* \geq c^*$ will follow from the next lemma, analogous to Lemma 6.5 of \cite{BH1}:

\begin{lem}
Let $(u,c)$ be a solution to problem (\ref{KP1}) with $u_s > 0$. Then there exists a $\lambda > 0$ and a function $\phi(y,\tau) > 0$ such that the function $w = e^{\lambda s}\phi(y,\tau)$ solves
\br
L_c w + f'(0) w = 0 .
\er 
\end{lem}
\nit {\it Proof:}
The proof adapts Lemma 6.5 of \cite{BH1}. Specifically, we define
\be
0 \leq \lambda = \liminf_{s \to -\infty} \frac{u_s(s,y,\tau)}{u(s,y,\tau)} ,
\ee
and we let $(s_j,y_j,\tau_j)$ be a sequence such that
\be
\lim_{j \to \infty} \frac{u_s(s_j,y_j,\tau_j)}{u(s_j,y_j,\tau_j)} = \lambda .
\ee
By the compactness of $T^n$, we may assume $(y_j,\tau_j) \to (y_\infty,\tau_\infty)$ for some point $(y_\infty,\tau_\infty) \in T^n$. As in \cite{BH1}, the family of functions
\be
w_j(s,y,\tau) = \frac{u(s + s_j,y,\tau)}{u(s_j,y_j,\tau_j)}
\ee
satisfies $w_j > 0$ and
\be
L_c w_j + \frac{f(u(s+s_j,y,\tau)}{u(s+s_j,y,\tau)} w_j = 0 .
\ee
Parabolic estimates imply that we can take a subsequence of the $w_j$ converging in $C^1_{\text{loc}}(R \times T^n)$ to a function $w \geq 0$ that solves
\be
L_c w + f'(0)w =  0 . \label{VAR11}
\ee
Also, $w(0,y_\infty,\tau_\infty) = 1$, so the maximum principle implies $w > 0$. Now we will see that $w$ has the form  $w = e^{\lambda s}\phi(y,\tau)$. Note that for each index $j$, 
\be
(w_j)_s(s,y,\tau) = \frac{u_s(s + s_j,y,\tau)}{u(s + s_j,y,\tau)} w_j(s,y,\tau) ,
\ee
so that $w_s(s,y,\tau) \geq \lambda w(s,y,\tau)$ with $w_s(0,y_\infty,\tau_\infty) = \lambda w(0,y_\infty,\tau_\infty) = \lambda$. Thus, the function defined by $z(s,y,\tau) = \frac{w_s}{w}$ satisfies $z \geq \lambda $ for all $(s,y,\tau) \in R \times T^n$ with equality holding at the point $(0,y_\infty,\tau_\infty)$. A straightforward computation shows that $z$ solves
\be
L_c z + 2 \left(\frac{\nabla_{s,y}w}{w}\cdot \nabla_{s,y}z \right) = 0 \quad \text{for} \quad (s,y,\tau) \in R \times T^n .
\ee
The strong maximum principle then implies that $z \equiv \lambda$, so that 
\br
\ds \left( w e^{-\lambda s} \right)  &=& w_s e^{-\lambda s} - \lambda w e^{-\lambda s} \no \\
&=& w \left( z e^{-\lambda s} - \lambda e^{-\lambda s} \right) \equiv 0 .
\er
Hence, the function $\phi = w(s,y,\tau)e^{- \lambda s} = \phi(y,\tau)$ is independent of $s$ and is positive. We conclude that $w(s,y,\tau) = e^{\lambda s}\phi(y,\tau)$. To show that $\lambda \neq 0$, we plug the function $w$ back into equation (\ref{VAR11}) and find that $\phi(y,\tau) > 0$ solves
\be
L^\lambda \phi = \Delta_y \phi + [\lambda^2 + \lambda b(y,\tau) + f'(0)]\phi - \phi_\tau = -c \lambda \phi .
\ee
From this equation we see that if $\lambda= 0$ then $\phi$ must be a constant, by the maximum principle, hence $w$ must be a constant. 
However, $w$ cannot be constant since $w(0,y_\infty,\tau_\infty) = 1$ and $w$ solves (\ref{VAR11}) with $f'(0) > 0$. This proves the lemma.

Applying this result to the solution $(u^*,c^*)$, we obtain a function $w = e^{\lambda s}\phi(y,\tau)$, where $\phi(y,\tau) > 0$ solves
\be
L^\lambda \phi = \Delta_y \phi + [\lambda^2 + \lambda b(y,\tau) + f'(0)]\phi - \phi_\tau = -c^* \lambda \phi ,
\ee
so that $\phi(y,\tau)$ is the principal eigenfunction of (\ref{VAR1}) with principal eigenvalue $\mu(\lambda) = -c^*\lambda$. Thus 
\be
-c^* = \frac{\mu(\lambda)}{\lambda} \geq \inf_{\lambda > 0}\frac{\mu(\lambda)}{\lambda} = -\gamma^* ,
\ee
so that $c^* \leq \gamma^*$. We have now shown that $c^* = \gamma^*$.

\begin{prop}
For $\lambda > 0$, the function $h(\lambda) = \frac{\mu(\lambda)}{\lambda}$ attains its infimum $-c^*$ at a unique point $\lambda^*>0$. Moreover, $h(\lambda)$ is decreasing in the region $(0,\lambda^*)$ and is 
increasing in the region $(\lambda^*,+\infty)$. Thus, no other local minima exist. 
\end{prop}
\nit {\it Proof:}
 As with (\ref{EX2}) and (\ref{EX3}), the function $h(\lambda)$ can be expressed as
\be
h(\lambda) = \frac{\mu(\lambda)}{\lambda} = \lam  + \frac{1}{\lambda}f'(0) + \frac{\int_{T^n} b\phi}{\int_{T^n}\phi} ,
\ee
where $\phi$ depends on $\lambda$. The last term is bounded by 
$\|b\|_{\infty}$ uniformly in $\lambda $ thanks to $\phi > 0$. 
From this we see that $h(\lambda) \to \infty$ as $\lambda \to 0$ 
and as $\lambda \to +\infty$.

Now we make use of some results from \cite{Hess}. Suppose that for some $c \in R$, $h(\lambda)=-c$ from some $\lambda$. Let $\phi > 0$ be an 
eigenfunction associated with principal eigenvalue $\mu(\lambda)$. 
Then, by (\ref{VAR1}) we have
\be
\Delta_y \phi + [\lambda^2 + \lambda (b(y,\tau) + c) + f'(0)]\phi - \phi_\tau = 0 \label{VAR12} 
\ee
since $c = -\frac{\mu(\lambda)}{\lambda}$. Therefore, $h(\lambda) = -c$ if and only if $\rho_c(\lambda) = -\lambda^2 - f'(0)$, where $\rho_c(\lambda)$ is the principal 
eigenvalue and $\phi>0$ the principal eigenfunction of the operator $\hat L$ defined by
\be
\hat L \phi = \Delta_y \phi + \lambda (b(y,\tau) + c)\phi - \phi_\tau = \rho_c(\lambda) \phi . \label{VAR13}
\ee
According to Lemma 15.2 of \cite{Hess}, the function $\lambda \mapsto \rho_c(\lambda)$ is convex. Since the function $-\lambda^2 - f'(0)$ is strictly concave, there can 
be at most two values of $\lambda$ such that $\rho_c(\lambda) = -\lambda^2 -f'(0)$. Hence, for any c, there can be at most two values of $\lambda$ such that $h(\lambda) = -c$.  This fact, 
combined with the continuity of $h(\lambda)$ and the limits $h(\lambda) \to \infty$ as $\lambda \to 0, \infty$, implies the Theorem.

The variational principle allows the comparison of $c^*$ 
with and without an additive time dependent component:
\begin{prop}
Write $b(y,\tau) = b_0 (y) + b_{1}(y,\tau)$, where 
$b_0(y) = \int_{0}^{1}\, b(y,\tau) \, d\tau $, 
$\int_{0}^{1}\, b_1(y,\tau) \, d\tau =0$. Let $c^{*}_{sp}$ be the 
minimal speed when $b_1 \equiv 0$, then $|c^{*}| \geq |c^{*}_{sp}|$, 
equality holds if $b_{1}=b_{1}(\tau)$. 
\end{prop}
This follows from the analogous comparison property of $\mu (\lambda)$, 
see Theorem 2.1 of \cite{HShenV}. 
 
\section{Variational Principle Based Computation}
\begin{figure}[tb]
\centerline{\epsfig{file=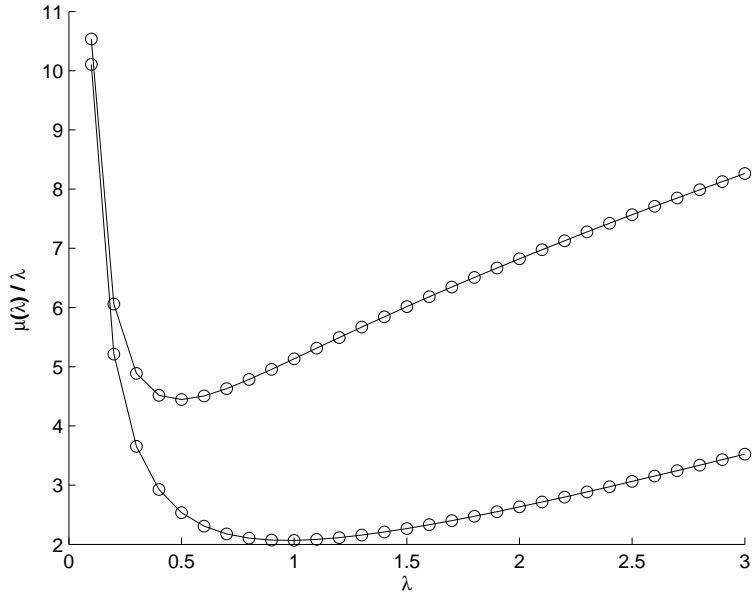,width=330pt}}
\caption{Two curves $\lambda \mapsto \frac{\mu(\lambda)}{\lambda}$}
\label{fig1}
\end{figure}

\setcounter{equation}{0}
The variational formula of Theorem 4.1 is useful from a 
computational point of view because it allows us to compute the 
minimal speed $c^*$ without computing the solution $u$ of (\ref{e1}) 
over long times. 
The variational formula, on the other hand, 
is a simple minimization problem in just one dimension: 
for a given shear flow, we want to compute 
\be
-c^* = \inf_{\lambda >0} \frac{\mu(\lambda)}{\lambda} . \no
\ee
By Lemma 4.2, the minimum is attained at a unique point $\lambda^*$. For a given $\lambda>0$, we find $\mu(\lambda)$ by solving the eigenvalue problem
\be
\Delta_y \phi + [\lambda^2 + \lambda b(y,\tau) + f'(0)]\phi = \mu(\lambda) \phi . \no
\ee
We approximate the infimum along the curve 
$\lambda \mapsto \frac{\mu(\lambda)}{\lambda}$ using a Newton's method 
with line search. In this way, 
our approximation converges quadratically in the region 
near the infimum and decreases with each iteration. 
Two illustrative curves 
$\lambda \mapsto \frac{\mu(\lambda)}{\lambda}$ 
are shown in Figure \ref{fig1} for two different shears. 

Because of the periodicity and smoothness of the coefficients, 
we compute the principal eigenvalues $\mu(\lambda)$ using a 
Fourier spectral method, as described in \cite{Tref}.  
Consider $n=2$, $y \in [0,2\pi]$ and $\tau \in [0,2\pi]$. 
For a domain $[0,2\pi]$ we divide 
the space into N equally spaced grid points and approximate
\be
u(y) \approx \sum_{j=1}^N \bar u_j S_N(y-y_j)
\ee
for $h = \frac{2\pi}{N}$ and $y_j = jh$. The function $S_N(y)$ is the band-limited interpolant of the Dirac mass $\delta$:
\be
S_N(y) = \frac{\sin(\pi y/h)}{(2\pi /h)\tan(y/2)} .
\ee
Differentiation then translates into matrix multiplication:
\be
(Du)(y_i) \approx D_N \bar u 
\ee
where $D_N$ is the matrix with entries
\br
(D_{N})_{ij} = \left\{ \begin{array}{ll} 0 \;\; &  \text{if}\;\;i=j \\ 
\frac{1}{2}(-1)^{(i-j)} \cot \left(\frac{(i-j)h}{2}\right)\;\; & \text{if}\;\; i \neq j , \end{array} \right.
\er
and for a second derivative
\br
(D^{(2)}_{N})_{ij} = \left\{ \begin{array}{ll} -\frac{\pi^2}{3h^2}-\frac{1}{6} \;\; &  \text{if}\;\;i=j \\ 
-\frac{(-1)^{(i-j)}}{2 \sin^2 ((i-j)h/2)} \;\; & \text{if}\;\; i \neq j . \end{array} \right.
\er
Performing a similar discretization in the $\tau$ dimension, we see that the numerical eigenvalue problem becomes
\be
A \bar u = (D^{(2)}_{N,y} + B_\lambda - D_{N,t}) \bar u= \mu(\lambda) \bar u ,
\ee
where $B_\lambda$ is the diagonal matrix given by
\be
B_\lambda = (\lambda^2 + f'(0))I + \lambda \, \text{diag}(b),
\ee
and $\text{diag}(b)$ denotes the diagonal matrix whose $i^{\text{th}}$ diagonal entry is $b(y(i),\tau(i))$, the value of the function at the point
\be
y(i) = ((i+1)\mod N) h \;,\quad \tau(i) = \left \lfloor \frac{i}{N} \right \rfloor h .
\ee
In all simulations, we take $f(u) = u(1-u)$. The resulting matrix A is not symmetric, since the matrix $D_N$ is not symmetric. We use general double precision LAPACK routines \cite{LAP} to find the princpal 
eigenvalue, corresponding to a positive eigenvector, and thus we obtain points on the curve $\lambda \mapsto \frac{\mu(\lambda)}{\lambda}$.

\subsection{Example 1}
First, we define
\be
b_\delta(y,\tau) = \delta \sin(2 \pi y)(1.0 + \sin(2 \pi \tau)), \label{Num1}
\ee
and we let $c^*(\delta)$ denote the corresponding minimal speed.  We consider how $c^*(\delta )$ scales with the shear amplitude $\delta$, as 
in \cite{NX1}.  It is known that for any $\delta > 0$, $c^*(\delta) < c^*$ (Recall that in our notation, $c^* < 0$). Supposing $c^*(\delta) \approx c^*(0) + O(\delta^p)$ we wish to numerically compute the exponent $p$.  According to \cite{NX1}, we 
expect $c^*(\delta ) - c^*_0$ is $O(\delta^2)$ for $\delta $ small, and $O(\delta)$ for $\delta $ large.  

\begin{figure}[tb]
\centerline{\epsfig{file=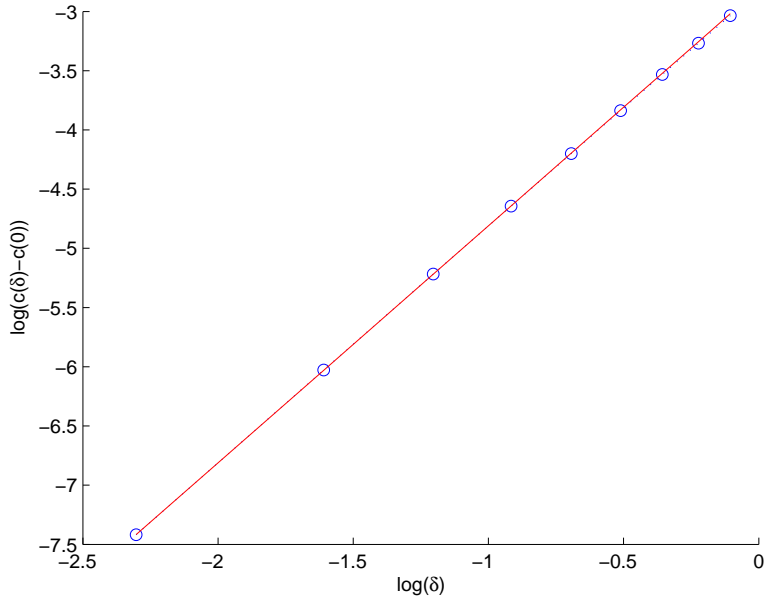,width=330pt}}
\caption{Log-Log plot of enhancement. Slope = 1.98 in the 
small shear amplitude regime.}
\label{fig3}
\end{figure}

Using the numerically calculated speeds $c^*(\delta)$, we determined the exponents $p$ using the least squares method to fit a line to a log-log plot 
of speed versus amplitude.  That is, we determined the slope of the best-fit line 
through the data points $(\log(\delta), \log(c^*(0)-c^*(\delta)))$ for 
each shear amplidute $\delta$.  Figure \ref{fig3} shows this log-log plot. We find that for small 
amplitude shears the minimal speed grows quadratically with the shear amplitude $\delta$: $p \approx 1.98$. For large amplitude 
shears, the enhancement is linear: $p \approx 1.05$, as 
shown next in Example 2. These results agree with 
\cite{NX1}, but are more accurate and computed much faster.

\subsection{Example 2}
\begin{figure}[tb]
\centerline{\epsfig{file=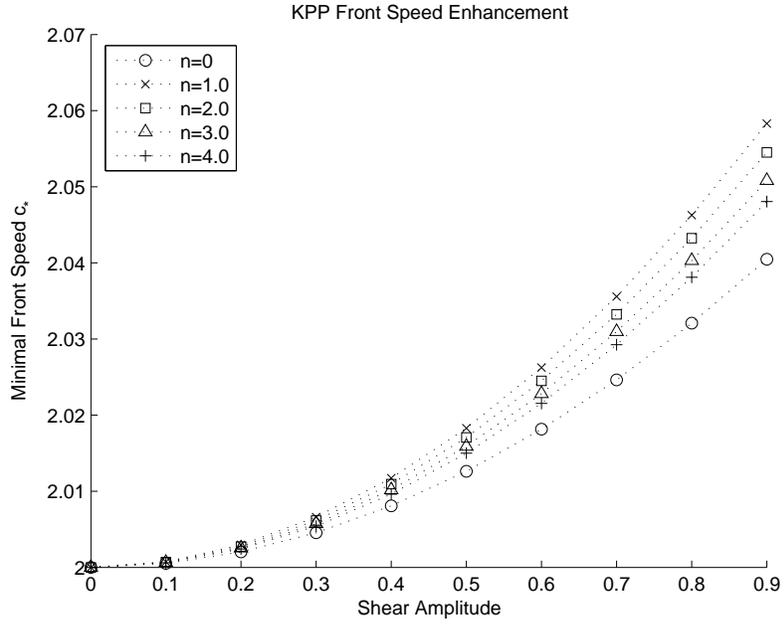,width=330pt}}
\caption{Enhancement by shears with various temporal frequencies,
quadratic speed growth in the small shear amplitude regime.}
\label{fig2}
\end{figure}

\begin{figure}[tb]
\centerline{\epsfig{file=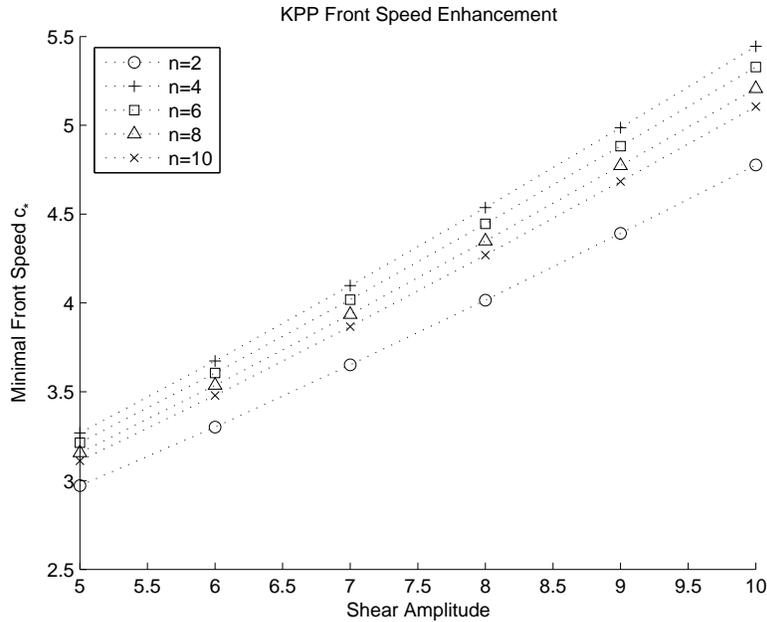,width=330pt}}
\caption{Enhancement by shears with various temporal frequencies, 
linear speed growth in the large shear amplitude regime.}
\label{fig4}
\end{figure}

Now we consider the effect of temporal frequency on the speed enhancement. For each $n \in \{ 0, 1, 2, 3, 4 \}$, we define
\be
b_{\delta,n}(y,\tau) = \delta \sin(2 \pi y)(1.0 + \sin(2 \pi n \tau)) . \label{Num2}
\ee
Figure \ref{fig2} shows the curves $\delta \mapsto -c^*(\delta,n)$ for each $n$, plotted simultaneously, for small amplitudes.  Figure \ref{fig4} shows the same plot for larger amplitudes. Notice that the curves become linear for $\delta$ sufficiently large. When $n=0$, the shear is independent of time, and the simulations show that the enhancement from $n=0$ is weakest (cf. Proposition 4.2). 
For time-dependent shears, the enhancment decreases monotonically with increasing temporal frequency.

Note that when $\delta = 0$, $b_{\delta,n}(y,t)\equiv 0$ and the minimal speed is $c^* = -2.0$, agreeing with the well-known formula
\be
c^* = - 2\sqrt{f'(0)} = -2. \no
\ee

\newpage
\section{Acknowledgements}
The work was partially supported by NSF grant ITR-0219004. 
J. X. would like to acknowledge the 
Faculty Research Assignment Award at the University of Texas at Austin, and
a fellowship from the John Simon Guggenheim Memorial Foundation.

\bibliographystyle{plain}

\begin{thebibliography}{99}
\bibitem{LAP} E. Anderson, et. al., 
{\em LAPACK Users' Guide}, Third Edition, Society for Industrial and Applied Mathematics, 1999.

\bibitem{BH1} H. Berestycki and F. Hamel,
{\em Front propagation in periodic excitable media}, Comm. in Pure and Applied Math. Vol. 60, (2002), pp. 949-1032.

\bibitem{BerN1} H. Berestycki, L. Nirenberg, 
{\em Travelling fronts in cylinders},  Ann. Inst. H. Poincare, 
Anal. Nonlineaire, Vol. 9, (1992), No. 5, pp. 497-572.

\bibitem{CW} P. Clavin, and F. A. Williams,
{\em Theory of premixed-flame propagation in large-scale
turbulence}, J. Fluid Mech., 90, (1979), pp. 598-604.

\bibitem{Const1} P. Constantin, A. Kiselev, A. Oberman, L. Ryzhik,
{\em Bulk burning rate in passive-reactive diffusion}, Arch Rat.
Mech Analy, 154, (2000), pp. 53-91.

\bibitem{HPS}S. Heinze, G. Papanicolaou, A. Stevens,
{\em Variational principles for propagation speeds in inhomogeneous media},
SIAM J. Applied Math, 62, no. 1, (2001), pp. 129 - 148.

\bibitem{Hess} P. Hess, 
{\em Periodic-Parabolic Boundary Value Problems and Positivity}, 
John Wiley and Sons, 1991. 

\bibitem{HShenV} V. Hutson, W. Shen and G. T. Vickers,
{\em Estimates for the principal spectrum point for certain 
time-dependent parabolic operators},
Proceedings of the American Mathematical Society, 129 (2001), pp. 1669-1679.

\bibitem{KA}A. R. Kerstein and W. T. Ashurst,
{\em Propagation rate of growing interfaces in stirred fluids},
Phys. Rev. Lett. 68, (1992), p. 934.

\bibitem{KR}A. Kiselev, L. Ryzhik,
{\em Enhancement of the traveling front speeds in reaction-diffusion
equations with advection}, Ann. de l'Inst. Henri Poincar\'e,
Analyse Nonlin\'eaire, 18, (2001), pp. 309--358.

\bibitem{Lad} O.A. Ladyzenskaja, V.A. Solonnikov, N.N. Ural'ceva, 
{\em Linear and Quasi-linear Equations of Parabolic Type}, 
American Mathematical Society, 1968. 

\bibitem{Laz} A.C. Lazer, {\em Some remarks on periodic solutions of parabolic differential equations} In Dynamical Systems II, A. Bednarek and L. Cesari, eds.  
Academic Press, New York, (1982), pp. 227-246.

\bibitem{Lbm} G. Lieberman, 
{\em Second Order Parabolic Differential Equations}. 
World Scientific, 1996. 

\bibitem{Kh}B. Khouider, A. Bourlioux, A. Majda,
{\em Parameterizing turbulent flame speed-Part I: unsteady shears, flame
residence time and bending}, Combustion Theory and Modeling, 5 (2001),
pp. 295-318.

\bibitem{MS1}A. Majda, P. Souganidis,
{\em Large scale front dynamics for turbulent reaction-diffusion equations
with separated velocity scales}, Nonlinearity, 7 (1994), pp. 1-30.

\bibitem{MS2}A. Majda and P. Souganidis,
{\em Flame fronts in a turbulent combustion model with fractal velocity
Fields}, Comm Pure Appl Math, Vol. LI (1998), pp. 1337-1348.

\bibitem{NX1}J. Nolen, J. Xin,
{\em Reaction diffusion front speeds in spacially-temporally periodic 
shear Flows}, SIAM J. Multiscale Modeling and Simulation, 
Vol. 1, (2003), No. 4, pp. 554-570.

\bibitem{PX}G. Papanicolaou, J. Xin,
{\em Reaction-diffusion fronts in periodically layered media},
J. Stat. Physics, 63 (1991), pp. 915-931.

\bibitem{PW}M. Protter and H. Weinburger, {\em Maximum Principles 
in Differential Equations}. Prentice Hall, 1967. 

\bibitem{Ro}P. Ronney,
{\em Some open issues in premixed turbulent combustion},
in: Modeling in Combustion Science (J. D. Buckmaster and T. Takeno, Eds.),
Lecture Notes In Physics, Vol. 449, Springer-Verlag, Berlin, (1995), pp. 3-22.

\bibitem{Tref}L. N. Trefethen,
{\em  Spectral Methods in MATLAB}, Society for Industrial and 
Applied Mathematics, 2000.

\bibitem{Vlad}N. Vladimirova, P. Constantin, A. Kiselev, O. Ruchayskiy, 
L. Ryzhik, 
{\em Flame enhancement and quenching in fluid flows}, Combust. Theory Model., 
7, (2003), pp. 487-508.

\bibitem{Xin0}J. Xin,
{\em Existence and stability of travelling waves in
periodic media governed by a bistable nonlinearity,} J. Dynamics
Diff. Eqs., 3, (1991), pp. 541--573.

\bibitem{Xin3}J. Xin,
{\em Existence of planar flame fronts in convective-diffusive
periodic media}, Arch. Rat. Mech. Anal., 121, (1992), pp. 205--233.

\bibitem{Xin1} J. Xin,
{\em Front propagation in heterogeneous media,}
SIAM Review, Vol. 42, No. 2, June 2000, pp. 161-230. 

\bibitem{Xin2} J. Xin,
{\em KPP front speeds in random shears and the parabolic Anderson problem}, 
Methods and Applications of Analysis, Vol. 10, No. 2, (2003), pp. 191-198.

\bibitem{Yak}V. Yakhot,
{\em Propagation velocity of premixed turbulent flames},
Comb. Sci. Tech 60, (1988), p. 191.

\end{thebibliography}

\end{document}